# The Mod 2 Cohomology Rings of Rank 3 Simple Groups are Cohen–Macaulay

Alejandro Adem[†]
Department of Mathematics
University of Wisconsin
Madison WI 53706

R. James Milgram[†]
Department of Mathematics
University of New Mexico
Albuquerque NM 87131

Dedicated to Professor William Browder

## §0 INTRODUCTION

In recent years we have been studying the cohomology rings of the sporadic simple groups and their relations with problems in algebra and homotopy theory. This work was motivated in part by results of D. Quillen [Q3] and J. Carlson [C], connecting $H^*(G; \mathbb{K})$ ($\mathbb{K}$ a field of characteristic $p$) with the structure of modular representations of $G$, and also by the connection of finite group theory with stable homotopy theory through the identification $B_{\mathcal{S}_\infty}^+ \simeq Q(S^0)$.

The rank (or 2-rank) of a finite group $G$ is the dimension of the largest 2-elementary subgroup $(\mathbb{Z}/2)^m$ contained in $G$, and this provides, at least in low rank, a convenient method for organizing the simple groups. The only simple group of rank 1 is $\mathbb{Z}/2$. For the simple groups of rank two the classification is given by the following theorem due to Alperin, Brauer, Gorenstein, Lyons, and Walter.

THEOREM: *The simple groups of rank 2 are $PSL_2(\mathbb{F}_q)$, $q$ odd and $q \geq 5$, $PSL_3(\mathbb{F}_q)$, $PSU_3(\mathbb{F}_q)$ for $q$ odd, $PSU_3(\mathbb{F}_4)$, the alternating group $\mathcal{A}_7$, and the Mathieu group, $M_{11}$.* ∎

For the simple groups of rank 3 the classification is given by O'Nan and Stroth.

THEOREM: *The simple groups of rank 3 are $G_2(q)$, $^3D_4(q)$, for $q$ odd, $^2G_2(3^n)$, $n$ odd, $n \geq 3$, $PSL_2(\mathbb{F}_8)$, $Sz(8)$, $PSU_3(\mathbb{F}_8)$, $J_1$, $M_{12}$, $O'N$.* ∎

(We review these classification results in §1.)

In this list, the mod(2) cohomology of the groups of Lie type *when the characteristic $q$ is odd* are closely tied to the cohomology of the classifying space $B_G$, where $G$ is the associated complex group, at least when the resulting graph automorphism allows a twisting over $\mathbb{C}$; this happens to be the case for $PSU_3$ but not for the other twisted groups in the list above (see [Q1], [Q2]). However, the cohomology of the triality twisted groups $^3D_4(q)$ is determined in [FM], and, since $Syl_2(^2G_2(3^n))$ is elementary abelian $H^*(^2G_2(3^n); \mathbb{F}_2) \cong H^*(Syl_2(^2G_2(3^n)); \mathbb{F}_2)^N$ where $N = \mathbb{Z}/7 \times_T \mathbb{Z}/3$ is the Weyl group of $(\mathbb{Z}/2)^3$ in $^2G_2(3^n)$. Note that this is also isomorphic to $H^*(J_1; \mathbb{F}_2)$ (see [W]).

[†] Both authors were partially supported by grants from the National Science Foundation and by the ETH–Zürich; the first author was also supported by an NSF Young Investigator Award.

The mod 2 cohomology of $\mathcal{A}_7$ is quite easy since $Syl_2(\mathcal{A}_7) = D_8$, the dihedral group; the result is given in several sources (see [AM1]). $H^*(M_{11}; \mathbb{F}_2)$ was originally determined by Webb [W] and Benson and Carlson [BC], though a much simpler derivation is given in [AM1]. Similarly, $H^*(M_{12}; \mathbb{F}_2)$ was determined in [AMM], and recently we determined $H^*(O'N; \mathbb{F}_2)$ ([AM4] and [M]).

Thus, the complete determination of the cohomology rings $H^*(G; \mathbb{F}_2)$ for all simple groups of 2-rank $\leq 3$ reduces to the study of the groups $PSU_3(\mathbb{F}_{2^n})$ and $Sz(2^n)$. Since we are working here at the characteristic we cannot take advantage of the étale techniques which were so successful previously. However, it turns out that the Sylow 2-subgroups in these cases are sufficiently simple that we can study their cohomology rings anyway. Although we are specifically concerned with $PSU_3(\mathbb{F}_4)$, $PSU_3(\mathbb{F}_8)$ and $Sz(8)$, we provide an analysis for the general situation.

The cohomology ring of a finite group is often quite complicated, which makes it exceedingly hard to obtain general structural results. Perhaps the nicest property one could hope for is that $H^*(G, \mathbb{F}_p)$ satisfy the *Cohen–Macaulay condition*: that is to say, it is free and finitely generated over a polynomial subalgebra of finite rank. In fact we will prove the following result.

THEOREM: *If $G$ is a finite simple group and its 2-rank is at most three, then $H^*(G, \mathbb{F}_2)$ is Cohen–Macaulay.* ∎

The theorem above is a corollary of the general result that is the culmination of the arguments in sections 2 and 3:

THEOREM: *If $G$ is one of the simple groups of Lie type, $Sz(2^{2n+1})$ or $PSU_3(\mathbb{F}_{2^n})$, then $H^*(G; \mathbb{F}_2)$ is Cohen-Macaulay.* ∎

In turn, this follows from the determination, in §2, of the structure of $Syl_2(G)$ for $G$ as above as a central extension, with center an elementary 2-group, and showing that for these $Syl_2(G)$, any element of order two is contained in the center. In §3 we study the cohomology rings of finite 2-groups with the property above and prove that $H^*(G; \mathbb{F}_2)$ is always Cohen-Macaulay, thus proving our result. We would like to point out that after this paper was written we learned that J. Clark [Cl] in his Oxford Ph. D. thesis has explicitly computed the mod 2 cohomology of $PSU_3(\mathbb{F}_4)$ and $Sz(8)$. Hence the mod 2 cohomology of practically all low rank simple groups has been explicitly calculated.

In the final section we describe how this result breaks down for groups of rank 4,

as illustrated for example by the Mathieu group $M_{22}$ or $PSL_3(\mathbb{F}_4)$. Work in progress indicates that a better grasp of the rank 4 situation will soon be available, but it will be of far greater complexity than what we present here.

This paper is dedicated to Bill Browder in celebration of his sixtieth birthday. In addition to his well–known geometric insight, he has had a long-standing interest in the cohomology of finite groups, which was a source of motivation for both authors. We owe him a debt of gratitude for his advice and encouragement.

## §1 LOW RANK SIMPLE GROUPS

In this section we will recall a few basic facts about simple groups. The main reference for them are Gorenstein's book [G1] and the Atlas [Co].

Let $p$ be a prime dividing the order of a finite group $G$; the $p$–rank of $G$ is by definition the maximal rank of a $p$–elementary abelian subgroup in $G$. As we will be dealing exclusively with the case $p = 2$, we will simply call this the rank of $G$, and abbreviate it by $r(G)$. Note that as a consequence of the Feit–Thompson Theorem, all simple groups are of *even order*, and hence $r(G) > 0$ if $G$ is simple.

To begin, we first observe that there are no simple groups of rank one except $\mathbb{Z}/2$. In this situation a 2–Sylow subgroup of $G$ is either cyclic or generalized quaternion. In the first case, Burnside's transfer theorem implies that $G$ has a normal 2–complement and hence is solvable, whereas in the other case, the Brauer–Suzuki theorem shows that $G$ has a unique, hence central, involution.

For the rank two case, a theorem of J. Alperin implies that $Syl_2(G)$ is either dihedral, quasi–dihedral, wreathed, or isomorphic to $Syl_2(PSU_3(\mathbb{F}_4))$. This was used as a basis for the classification of simple groups of rank 2, which we now describe.

THEOREM 1.1: (Gorenstein–Walter) *If $G$ is a finite simple group with a dihedral 2–Sylow subgroup, then it is isomorphic to one of the following groups:*

$PSL_2(\mathbb{F}_q)$, $q > 3$ *odd or* $A_7$. ∎

Next we have

THEOREM 1.2: (Alperin–Brauer–Gorenstein) *If $G$ is a finite simple group with a quasi–dihedral or wreathed 2–Sylow subgroup, then it is isomorphic to one of the following groups:*

$PSL_3(\mathbb{F}_q)$, $PSU_3(\mathbb{F}_q)$, $q$ *odd or* $M_{11}$, *the first Mathieu group.* ∎

Finally, we have the result from Lyons' thesis

THEOREM 1.3: (Lyons) If $G$ is a simple group with $Syl_2(G) \cong Syl_2(PSU_3(\mathbb{F}_4))$, then $G \cong PSU_3(\mathbb{F}_4)$. ∎

Combining these results, we obtain the following

(1.4) CLASSIFICATION OF RANK 2 SIMPLE GROUPS:

If $G$ is a simple group of rank 2, then $G \cong PSL_2(\mathbb{F}_q)$, ($q$ odd, $q \geq 5$), $PSL_3(\mathbb{F}_q)$ or $PSU_3(\mathbb{F}_q)$ ($q$ odd), $PSU_3(\mathbb{F}_4)$, $A_7$ or $M_{11}$. ∎

For the rank 3 case there is a similar result, proved by O'Nan and Stroth

(1.5) CLASSIFICATION OF RANK 3 SIMPLE GROUPS :

If $G$ is a simple group of rank 3, then $G \cong G_2(q)$ or $^3D_4(q)$ (triality twisted) for $q$ odd, $^2G_2(3^n)$, $n$ odd, $n > 1$, $PSL_2(\mathbb{F}_8)$, $Sz(8)$, $PSU_3(\mathbb{F}_8)$, $J_1$, $M_{12}$ or $O'N$. ∎

Note that in the list above, the last three are sporadic groups, the rest are finite Chevalley groups, or twisted versions of them.

Now we briefly review the preceding classification results from the point of view of group cohomology. First we recall that the cohomology ring $H^*(G, \mathbb{F}_2)$ is said to be *Cohen–Macaulay* (abbreviated from now on as $CM$) if there exists a polynomial subalgebra $\mathcal{R} \subset H^*(G, \mathbb{F}_2)$ over which the ring is *free and finitely generated*. It follows from the finite generation in group cohomology that if $H^*(Syl_2(G), \mathbb{F}_2)$ is $CM$, then so is $H^*(G, \mathbb{F}_2)$. Also note that if $\mathbb{K}$ is an extension of $\mathbb{F}_2$, then $H^*(G, \mathbb{K}) \cong H^*(G, \mathbb{F}_2) \otimes \mathbb{K}$, so that the $CM$ property over $\mathbb{K}$ is equivalent to having it over the ground field.

We begin our analysis by first pointing out that the cohomology of the dihedral groups and of the "wreathed groups" $(\mathbb{Z}/2^n) \wr \mathbb{Z}/2$ is well-known (see [AM1], Chapters III and IV) and easily verified to be $CM$. Next we observe that the groups $PSL_3(\mathbb{F}_q)$ and $PSU_3(\mathbb{F}_q)$ ($q$ odd) are quotients of $SL_3(\mathbb{F}_q)$ and $SU_3(\mathbb{F}_q)$ (respectively) by subgroups of *odd order*, and hence have the same mod 2 cohomology. They are, consequently, $CM$. (We can also check this from the explicit calculations in [Q1], [FP].) Since $Syl_2(\mathcal{A}_7) = D_8$, it too is $CM$.

The ring $H^*(Syl_2(M_{11}); \mathbb{F}_2)$ is *not* $CM$. However, $H^*(M_{11}; \mathbb{F}_2)$ is explicitly determined in [AM1] and is $CM$. We conclude that in the rank two case, the only possible group not having a $CM$ cohomology ring is $PSU_3(\mathbb{F}_4)$.

We defer this case for the moment and turn to the rank 3 case: first note that the groups $^2G_2(3^n)$, $n > 1$ odd, $PSL_2(\mathbb{F}_8)$ and $J_1$ have $(\mathbb{Z}/2)^3$ as their 2–Sylow subgroup. Hence they have $CM$ cohomology rings. The Chevalley groups of finite type have been extensively studied using étale are well as other methods; in particular the explicit calculations in [Q2], [K] and [FM] show that $G_2(q)$ and $^3D_4(q)$ ($q$ odd) have $CM$ cohomology rings. In [AMM] it was shown that $H^*(M_{12}, \mathbb{F}_2)$ is $CM$, and very recently we verified by an explicit calculation [AM4] that the same is true for the O'Nan group O'N. Hence we conclude that the only two possible rank 3 simple groups with cohomology not satisfying $CM$ are $PSU_3(\mathbb{F}_8)$ and $Sz(8)$.

## §2 THE 2–SYLOW SUBGROUPS OF $PSU_3(\mathbb{F}_{2^n})$ and $Sz(2^{2n+1})$

In this section we study the structure of the groups $Syl_2(G)$ when $G = PSU_3(\mathbb{F}_{2^n})$ or $Sz(2^{2n+1})$. Our main objective will be to prove that in each case every element of order two in the Sylow subgroup is contained in its center.

**The groups $PSU_3(\mathbb{F}_{2^n})$**

Consider the non-singular Hermitian form over $\mathbb{F}_{2^{2n}}$

$$H = \begin{pmatrix} 0 & 1 & 0 \\ 1 & 0 & 0 \\ 0 & 0 & 1 \end{pmatrix}.$$

$U_3(\mathbb{F}_{2^n})$ is the subgroup of $GL_3(\mathbb{F}_{2^{2n}})$ consisting of matrices $A$ so that $AH\bar{A}^t = H$, where $(\bar{A}^t)_{i,j} = (A_{j,i})^{2^n}$. Then, for $A$ diagonal, $A = \begin{pmatrix} a & 0 & 0 \\ 0 & b & 0 \\ 0 & 0 & c \end{pmatrix}$, we see that $A \in U_3(\mathbb{F}_{2^n})$ if and only if $c\bar{c} = 1$, and $\bar{b} = a^{-1}$. Now, suppose that $A$ is central so $a = b = c$. Then $a^{2^n} = a^{-1}$ so $a^{2^n+1} = 1$, and conversely, if $a^{2^n+1} = 1$ then the diagonal matrix with $a$'s along the diagonal is contained in $U_3(\mathbb{F}_{2^n})$. On the other hand, if we concentrate on $SU_3(\mathbb{F}_{2^n})$, we must also have $a^3 = 1$, so $a$ is a third root of unity, and, since $a^{2^n} = a$ for $a \in \mathbb{F}_{2^n}$ we must have $n$ odd in order that $SU_3(\mathbb{F}_{2^n})$ contain a non-trivial center. Hence the simple group

$$PSU_3(\mathbb{F}_{2^n}) = \begin{cases} SU_3(\mathbb{F}_{2^n}) & \text{if } n \text{ is even,} \\ SU_3(\mathbb{F}_{2^n})/(\mathbb{Z}/3) & \text{if } n \text{ is odd.} \end{cases}$$

Now we construct the Sylow 2-subgroup of $PSU_3(\mathbb{F}_{2^n})$ and its normalizer explicitly, as this is needed to understand its cohomology.

For $\theta, \gamma \in \mathbb{F}_{2^{2n}}$, define the element $\mathcal{G}_{\theta,\gamma} \in SL_3(\mathbb{F}_{2^{2n}})$ as the matrix

$$\mathcal{G}_{\theta,\gamma} = \begin{pmatrix} 1 & \gamma & \bar{\theta} \\ 0 & 1 & 0 \\ 0 & \theta & 1 \end{pmatrix}.$$

Note that

$$\mathcal{G}_{\theta,\gamma} H \bar{\mathcal{G}}^t_{\theta,\gamma} = \begin{pmatrix} \gamma & 1 & \bar{\theta} \\ 1 & 0 & 0 \\ \theta & 0 & 1 \end{pmatrix} \begin{pmatrix} 1 & 0 & 0 \\ \bar{\gamma} & 1 & \bar{\theta} \\ \theta & 0 & 1 \end{pmatrix}$$

$$= \begin{pmatrix} \gamma + \bar{\gamma} + \theta\bar{\theta} & 1 & 0 \\ 1 & 0 & 0 \\ 0 & 0 & 1 \end{pmatrix},$$

so $\theta_{\theta,\gamma} \in SU_3(\mathbb{F}_{2^{2n}})$ if and only if $\gamma + \bar{\gamma} + \theta\bar{\theta} = 0$. However, given any $\theta \in \mathbb{F}_{2^{2n}}$ there is at least one $\gamma$ so that $tr(\gamma) = N(\theta) \in \mathbb{F}_{2^n}$. Any two such $\gamma$'s differ by an element in $\mathbb{F}_{2^n}$, so, in fact, there are exactly $2^n$ suitable choices for $\gamma$. We also have

$$\begin{pmatrix} 1 & \gamma & \bar{\theta} \\ 0 & 1 & 0 \\ 0 & \theta & 1 \end{pmatrix} \begin{pmatrix} 1 & \tau & \bar{\mu} \\ 0 & 1 & 0 \\ 0 & \mu & 1 \end{pmatrix} = \begin{pmatrix} 1 & \gamma + \tau + \bar{\theta}\mu & \bar{\mu} + \bar{\theta} \\ 0 & 1 & 0 \\ 0 & \mu + \theta & 1 \end{pmatrix}$$

so $\mathcal{G}_{\theta,\gamma} \mathcal{G}_{\mu,\tau} = \mathcal{G}_{\theta+\mu, \gamma+\tau+\bar{\theta}\mu}$. It follows that $K = \{\mathcal{G}_{\theta,\gamma} \mid \gamma + \bar{\gamma} + \theta\bar{\theta} = 0\}$ is a subgroup of $SU_3(\mathbb{F}_{2^n})$.

Define $Z \subset K$ to be the subgroup of $Z = \{\mathcal{G}_{0,\gamma} \mid \gamma \in \mathbb{F}_{2^n}\}$. Then $Z$ is central in $K$ with quotient isomorphic to $\mathbb{F}_{2^{2n}}^+$. Moreover, $Z \cong \mathbb{F}_{2^n}^+$, and we have that $|K| = 2^{3n}$, and $K$ is given as a central extension of the form

2.1 $$1 \longrightarrow \mathbb{F}_{2^n}^+ \overset{\triangleleft}{\longrightarrow} K \overset{\pi}{\longrightarrow} \mathbb{F}_{2^{2n}}^+ \longrightarrow 1.$$

The extension data for this extension is given by

$$\mathcal{G}_{\theta,\gamma}^2 = \mathcal{G}_{0,\theta\bar{\theta}}$$
$$[\mathcal{G}_{\theta,\gamma}, \mathcal{G}_{\mu,\tau}] = \mathcal{G}_{0,\bar{\mu}\theta+\mu\bar{\theta}}$$

and gives $K$ as an explicit central extension. Moreover, from the expression for the order of $PSU_n(\mathbb{F}_q)$ given in [Co], page $x$, we see that $K = Syl_2(PSU_3(\mathbb{F}_{2^n}))$.

We can make the extension data more explicit. Let $\rho \in \mathbb{F}_{2^{2n}}$ be a primitive $(2^n+1)^{st}$ root of unity. Then $\rho\bar{\rho} = 1$ and, if $\lambda = \rho + \bar{\rho}$ we have that $x^2 + \lambda x + 1$ is irreducible over $\mathbb{F}_{2^n}$ with roots $\rho$ and $\bar{\rho}$. Moreover, $\lambda$ cannot lie in any proper subfield of $\mathbb{F}_{2^l} \subset \mathbb{F}_{2^n}$ since, if it did, then $(2^n+1)$ would have to divide $2^{2l}-1$, and since $l$ divides $n$ this is impossible.

It follows that a basis for $\mathbb{F}_{2^n}$ over $\mathbb{F}_2$ is

$$\{1, \lambda, \lambda^2, \ldots, \lambda^{n-1}\},$$

and this completes to a basis for $\mathbb{F}_{2^{2n}}$ over $\mathbb{F}_2$ of the form

$$\{1, \lambda, \lambda^2, \ldots, \lambda^{n-1}, \rho, \lambda\rho, \lambda^2\rho, \ldots, \lambda^{n-1}\rho\}.$$

Specializing the relations above to this basis we have

$$[\mathcal{G}_{\lambda^i \rho,-}, \mathcal{G}_{\lambda^j,-}] = \mathcal{G}_{0,\lambda^{i+j+1}}$$
$$[\mathcal{G}_{\lambda^i \rho,-}, \mathcal{G}_{\lambda^j \rho,-}] = 1$$
$$[\mathcal{G}_{\lambda^i,-}, \mathcal{G}_{\lambda^j,-}] = 1$$

while

$$(\mathcal{G}_{r\lambda^i,-})^2 = \mathcal{G}_{0,\lambda^{2i}}$$
$$(\mathcal{G}_{\lambda^i,-})^2 = \mathcal{G}_{0,\lambda^{2i}},$$

which gives an explicit presentation of $K$.

As a direct consequence we have

LEMMA 2.2: $Z \subset K$ is the unique maximal elementary 2-subgroup of $K$. In particular, every element of order two in $K$ is contained in $Z$. ∎

We now turn our attention to the remaining case, the Suzuki groups.

**The group $Syl_2(Sz(2^{2n+1}))$**

From [G2], p.153, we see that $Syl_2(Sz(2^{2n+1}))$ is the subgroup of $GL_3(\mathbb{F}_{2^{2n+1}})$ with generators

$$\begin{pmatrix} 1 & 0 & 0 \\ a^\theta & 1 & 0 \\ b & a & 1 \end{pmatrix}$$

where $a,b \in \mathbb{F}_{2^{2n+1}}$, and $\theta = 2^s$ with $2^{2s} \cong 2 \bmod (2^{2n+1}-1)$. Thus $\theta = 4$ for $2^{2n+1} = 8$, the case of most interest to us, but more generally, $\theta$ can be taken to be $2^{n+1}$. Consequently $|Syl_2(Sz(2^{2n+1}))| = 2^{4n+2}$, and it is given as a central extension

$$1 \longrightarrow \mathbb{F}^+_{2^{2n+1}} \overset{\triangleleft}{\longrightarrow} Syl_2(Sz(2^{2n+1})) \overset{\pi}{\longrightarrow} \mathbb{F}^+_{2^{2n+1}} \longrightarrow 1.$$

Here, the extension data is

$$\{a\}^2 = a^{1+\theta}$$
$$[\{a\},\{b\}] = a^\theta b + ab^\theta$$

where $\{a\} = \begin{pmatrix} 1 & 0 & 0 \\ a^\theta & 1 & 0 \\ 0 & a & 1 \end{pmatrix}$, while $b = \begin{pmatrix} 1 & 0 & 0 \\ 0 & 1 & 0 \\ b & 0 & 0 \end{pmatrix}$ give explicit representatives for the two parts of the central extension. In particular, the squaring relation above shows that once more we have only one maximal elementary 2-subgroup in $Syl_2(Sz(2^{2n+1}))$, since, if $m = (2^{2n+1}-1, 1+2^{n+1})$, then $m$ divides both $2^{2n+1}-1$ and $2^{2n+2}-1$, so $m$ divides $2^{2n+2} - 2^{2n-1} = 2^{2n-1}$ and $m$ is 1. We have proved

LEMMA 2.3: The center of $Syl_2(Sz(2^{2n+1}))$ is an elementary 2-group, $2^{2n+1}$, and every element of order two in $Syl_2(Sz(2^{2n+1}))$ is contained in the center. ∎

## §3 THE COHOMOLOGY CALCULATIONS

We consider a finite 2-group $G$ with center $Z$, and suppose that $Z$ is 2-elementary for simplicity in what follows. Specifically, assume $Z = 2^n$. Multiplication induces a *homomorphism*
$$Z \times G \xrightarrow{\mu} G$$
and there is a finite $l \geq 1$ so that, for each $x \in H^*(Z; \mathbb{F}_2)$, $x^{2^l} \in res_G^Z(H^*(G; \mathbb{F}_2))$. In particular, the polynomial algebra $\tilde{\mathcal{P}} = \mathbb{F}_2[x_1^{2^l}, x_2^{2^l}, \ldots, x_n^{2^l}]$ is in the image of restriction from $H^*(G; \mathbb{F}_2)$, and it follows that there is a polynomial subalgebra $\mathcal{P} \subset H^*(G; \mathbb{F}_2)$ which injects isomorphically to $\tilde{\mathcal{P}}$ under restriction.

LEMMA 3.1: *Under the assumptions above, $H^*(G; \mathbb{F}_2)$ is free as a module over $\mathcal{P}$.*

PROOF: The composition of homomorphisms
$$G \xrightarrow{i_2} Z \times G \xrightarrow{\mu} G$$
is the identity, where $i_2$ is injection to the second factor. Thus $H^*(G; \mathbb{F}_2)$ is a direct summand of $H^*(Z \times G; \mathbb{F}_2)$ as a $\mathcal{P}$-module. Moreover, $H^*(Z \times G; \mathbb{F}_2)$ is free over $\mathcal{P}$ on generators
$$\left(H^*(Z; \mathbb{F}_2)/I(\tilde{\mathcal{P}})\right) \otimes H^*(G; \mathbb{F}_2).$$
Thus $H^*(G; \mathbb{F}_2)$ is projective over $\mathcal{P}$, but an easy induction on dimension then shows it is free as well. ∎

COROLLARY 3.2: *For each $n \geq 1$, the rings*
$$H^*(Syl_2(PSU_3(\mathbb{F}_{2^n})); \mathbb{F}_2), \qquad H^*(Syl_2(Sz(2^{2n+1})); \mathbb{F}_2)$$
*are $CM$. Consequently, the rings $H^*(PSU_3(\mathbb{F}_{2^n}); \mathbb{F}_2)$ and $H^*(Sz(2^{2n+1}); \mathbb{F}_2)$ are $CM$ as well.*

Proof: Let $G$ be one of the Sylow subgroups above. It suffices to show that $H^*(G; \mathbb{F}_2)$ is $CM$. We know that $H^*(G; \mathbb{F}_2)$ is finitely generated, and from the facts that an element is nilpotent if it does not restrict non-trivially to $H^*(Z; \mathbb{F}_2)$, and the restriction of $x^{2^l}$ is in the polynomial algebra $\tilde{\mathcal{P}}$, we see that at most a finite number of the generators and their products are needed to describe $H^*(G; \mathbb{F}_2)$ as a module over $\mathcal{P}$. ∎

REMARK 3.3: In the three cases of most interest to us here, $PSU_3(\mathbb{F}_4)$, $PSU_3(\mathbb{F}_8)$, and $Sz(8)$, it is possible to make the calculations explicit and from this obtain, at least in the case of $Sz(8)$, the actual cohomology rings (see [Cl]).

REMARK 3.4: The results in this section are actually a special case of a general theorem proved by J. Duflot [D]: *if $E \subset G$ is a maximal p–elementary abelian subgroup of $G$, then $H^*(C_G(E), \mathbb{F}_p)$ is Cohen–Macaulay.*

## §4 CONCLUSIONS AND FINAL REMARKS

Combining the results in the previous three sections, we deduce the following

THEOREM 4.1: *The mod 2 cohomology ring of every simple group of rank $\leq 3$ is Cohen–Macaulay.* ∎

This result would seem to indicate to the casual reader that perhaps simple groups have cohomology rings which are particularly nice. Unfortunately, this naive suggestion collapses at the next rank. Note that for $H^*(G, \mathbb{F}_2)$ to be $CM$ it is necessary that its maximal elementary abelian subgroups (at p=2) be of the same rank. This fails to be true for $M_{22}$, the next Mathieu group; it has maximal elementaries of ranks three and four. The cohomology of this group is quite complicated, but nevertheless computable (see [AM2]). On the other hand, $PSL_3(\mathbb{F}_4)$ is a simple group with maximal elementaries all of rank four, but its cohomology is not $CM$, as can be verified from the explicit calculation in [AM3].

The main result in this paper shows that low rank simple groups are particularly accessible, and from a topological point of view their classifying spaces can often be modelled (2–locally) using fibrations

$$F \longrightarrow BG \longrightarrow X$$

where the cohomology of $X$ embeds as a polynomial subalgebra realizing the $CM$ property and $F$ has the homotopy type of a finite complex. Examples of this can be found in [AM], page 279–80 and in [BW].

In contrast, the situation at rank four and beyond is a lot more complex, and we are only starting to obtain some idea (based on very recent calculations or work in progress) as to how involved the cohomology of simple groups can truly be.